\documentclass[11pt]{article}
\usepackage[dvipsnames]{xcolor}
\usepackage{amsfonts}
\usepackage{setspace}
\usepackage{bm}
\usepackage{bbm}
\usepackage{amssymb}
\usepackage{indentfirst}
\usepackage{geometry}
\usepackage{mathtools}

\DeclarePairedDelimiter\abs{\lvert}{\rvert}%
\DeclarePairedDelimiter{\floor}{\lfloor}{\rfloor}

\DeclareMathOperator{\HP}{HP}

\newcommand*{\citena}[1]{%
\begingroup
\color{Green}
\romannumeral-`\x 
\setcitestyle{numbers}%
\citep{#1}
\endgroup
\ignorespacesafterend
}

\newcommand*{\citesup}[1]{%
\begingroup
\color{Green}
\citep{#1}
\endgroup
\ignorespacesafterend
}

\newcommand*{\eqrefe}[1]{%
\begingroup
(\color{BrickRed}
\romannumeral-`\x 
\setcitestyle{numbers}%
\ref{eq:#1}%
\endgroup
)\ignorespacesafterend
}

\newcommand*{\secrefe}[1]{%
\begingroup
(\color{Aquamarine}
\romannumeral-`\x 
\setcitestyle{numbers}%
\ref{#1}%
\endgroup
)\ignorespacesafterend
}

\newtagform{Tags}[\textcolor{BrickRed}]{\color{Black}(}{)}

\usepackage[numbers, sort]{natbib}
\setcitestyle{super}

\newcommand{\ii}{\bm{i}}

\begin{document}
\title{Lerch's $\Phi$ and the Polylogarithm at the Positive Integers}
\date{June 15, 2020}
\author{Jose Risomar Sousa}
\maketitle
\usetagform{Tags}

\begin{abstract}
We review the closed forms of the partial Fourier sums associated with $\HP_k(n)$ from a previous paper and create an asymptotic expression for $\HP(n)$ as a way to obtain formulae for the full Fourier series (if $|b|<1$, one obtains a surprising pattern, $\HP(n) \sim H(n)-\sum_{k\ge 2}(-1)^k\zeta(k)b^{k-1}$). Finally, the derived Fourier series formulae are used to obtain a formula for the Lerch transcendent function, $\Phi(e^z,k,b)$, and by extension the polylogarithm, $\mathrm{Li}_{k}(e^{z})$, at the positive integers $k$.
\end{abstract}

\usetagform{Tags}

\tableofcontents

\section{Introduction}
Since the Basel problem in 1650, scholars have been eager to find closed forms for similar infinite series, especially Dirichlet series. In this manuscript, formulae for the Lerch transcendent function, $\Phi(e^z,k,b)$, and the polylogarithm, $\mathrm{Li}_{k}(e^{z})$, are created that hold at the positive integers $k$. Conversely, a formula for the Hurwitz zeta function at the negative integers, $\zeta(-k,b)$, is also created, to complement a formula at the positive integers produced in \citena{Hurwitz}.\\

The advantage of formulae that only hold at the positive integers is that they are likely to be simpler and easier to work with. It is an obvious statement if, for example, one thinks about the closed forms of the zeta function at the positive even integers, $\zeta(2k)$, and its general integral, valid for $\Re{(k)}>1$.\\

The formulae derived herein are based on new expressions created for the generalized harmonic progressions: 
\begin{equation} \nonumber
\HP_{k}(n)=\sum_{j=1}^{n}\frac{1}{(a\,n+b)^{k}} \text{,}
\end{equation}
\noindent which have been extensively studied in two previous papers, and vary depending on whether the parameters, $a$ and $b$, are integer\citesup{GHP} or complex\citesup{CHP}. When $a=1$ and $b=0$, one has a notable particular case, the generalized harmonic numbers, $H_k(n)$.\\

In \citena{GHP}, expressions were derived for the partial Fourier sums, $C^z_{k}(a,b,n)$ and $S^z_{k}(a,b,n)$, associated with $\HP_{k}(n)$, which are reproduced again in the next section, with a short description.\\

The objective of this paper is to obtain the limit of those expressions as $n$ approaches infinity, and then combine them to obtain the proper Lerch transcendent function, $\Phi$, at the positive integers.\\

The approach requires evaluating the limit of $\HP(n)-H(n)$, with $2b$ a non-integer complex number. Since this limit can also be attained by means of the digamma function, $\psi(n)$, this is just a new, more interesting way of deriving that limit.\\

In section \secrefe{int_lim}, the limits of the integrals that appear in the expressions of $\HP_k(n)$ as $n$ tends to infinity are reviewed, since they are central to this solution.\\

The process of obtaining the limits of $C^z_{2k}(b,n)$ and $S^z_{2k+1}(b,n)$ is much simpler than that of $C^z_{2k+1}(b,n)$ and $S^z_{2k}(b,n)$, since the latter involve the limit of $\HP(n)$, which diverges.

\section{The partial Fourier sums}
The subsequent expressions are the partial sums of the Fourier series associated with the generalized harmonic progressions from \citena{GHP}, and hold for all complex $z$, $a$ and $b$ and for all integer $n\geq 1$.\\

\indent By definition $\HP_{0}(n)=0$ for all positive integer $n$, so they actually have no effect in the sums. If $b=0$, any term that has a null denominator can be disregarded and the equation still holds (technically, the limit as $b$ tends to zero is taken, as seen in section \secrefe{part}). 

\subsection{$C^z_{2k}(a,b,n)$ and $S^z_{2k+1}(a,b,n)$} \label{Final_1}
For all integer $k \geq 1$:
\begin{multline} \label{eq:C^z_2k_parc}
\sum_{j=1}^{n}\frac{1}{(a j+b)^{2k}}\cos{\frac{2\pi(a j+b)}{z}}=-\frac{1}{2b^{2k}}\left(\cos{\frac{2\pi b}{z}}-\sum_{j=0}^{k}\frac{(-1)^j}{(2j)!}\left(\frac{2\pi b}{z}\right)^{2j}\right)
\\+\frac{1}{2(a n+b)^{2k}}\left(\cos{\frac{2\pi(a n+b)}{z}}-\sum_{j=0}^{k}\frac{(-1)^j}{(2j)!}\left(\frac{2\pi(a n+b)}{z}\right)^{2j}\right)
\\+\sum_{j=1}^{k}\frac{(-1)^{k-j}}{(2k-2j)!}\left(\frac{2\pi}{z}\right)^{2k-2j}\HP_{2j}(n)
\\+\frac{(-1)^k}{2(2k-1)!}\left(\frac{2\pi}{z}\right)^{2k}\int_{0}^{1}(1-u)^{2k-1}\left(\sin{\frac{2\pi(a n+b)u}{z}}-\sin{\frac{2\pi b u}{z}}\right)\cot{\frac{\pi a u}{z}}\,du 
\end{multline}
\indent For all integer $k \geq 0$:
\begin{multline} \label{eq:S^z_2k+1_parc}
\sum_{j=1}^{n}\frac{1}{(a j+b)^{2k+1}}\sin{\frac{2\pi(a j+b)}{z}}=-\frac{1}{2b^{2k+1}}\left(\sin{\frac{2\pi b}{z}}-\sum_{j=0}^{k}\frac{(-1)^j}{(2j+1)!}\left(\frac{2\pi b}{z}\right)^{2j+1}\right)
\\+\frac{1}{2(a n+b)^{2k+1}}\left(\sin{\frac{2\pi(a n+b)}{z}}-\sum_{j=0}^{k}\frac{(-1)^j}{(2j+1)!}\left(\frac{2\pi(a n+b)}{z}\right)^{2j+1}\right)
\\+\sum_{j=1}^{k}\frac{(-1)^{k-j}}{(2k+1-2j)!}\left(\frac{2\pi}{z}\right)^{2k+1-2j}\HP_{2j}(n)
\\+\frac{(-1)^k}{2(2k)!}\left(\frac{2\pi}{z}\right)^{2k+1}\int_{0}^{1}(1-u)^{2k}\left(\sin{\frac{2\pi(a n+b)u}{z}}-\sin{\frac{2\pi b u}{z}}\right)\cot{\frac{\pi a u}{z}}\,du 
\end{multline}

\subsubsection{The limits of $C^z_{2k}(n)$ and $S^z_{2k+1}(n)$} \label{facil}
For comparison purposes, let us review some limits that were derived previously for the particular cases $C^z_{2k}(n)$ and $S^z_{2k+1}(n)$ (that is, $a=1$ and $b=0$). The limits of the more general expressions are expected to coincide with them.\\

At infinity, these particular cases become Fourier series (denoted here by $C^z_{2k}$ and $S^z_{2k+1}$, respectively), which for positive or non-negative integer $k$ have limits given by:
\begin{equation} \label{eq:C^z_2k}
\sum_{j=1}^{\infty}\frac{1}{j^{2k}}\cos{\frac{2\pi j}{z}}=\sum_{j=0}^{k}\frac{(-1)^{k-j}}{(2k-2j)!}\left(\frac{2\pi}{z}\right)^{2k-2j}\zeta(2j)+\frac{(-1)^k\abs{z}}{4(2k-1)!}\left(\frac{2\pi}{z}\right)^{2k} 
\end{equation}
\begin{equation} \label{eq:S^z_2k+1}
\sum_{j=1}^{\infty}\frac{1}{j^{2k+1}}\sin{\frac{2\pi j}{z}}\\=\sum_{j=0}^{k}\frac{(-1)^{k-j}}{(2k+1-2j)!}\left(\frac{2\pi}{z}\right)^{2k+1-2j}\zeta(2j)+\frac{(-1)^k\abs{z}}{4(2k)!}\left(\frac{2\pi}{z}\right)^{2k+1}  
\end{equation}\\
\indent These limits hold for real $\abs{z}\ge 1$ ($k=0$ and $\abs{z}=1$ are exceptions and also trivial cases). For $S^1_{1}=0$ the formula breaks down (see section \secrefe{int_lim} to know why). Both of these results are known in the literature, they are rewrites of equations that feature in \citena{dois} (page 805).

\subsection{$C^z_{2k+1}(a,b,n)$ and $S^z_{2k}(a,b,n)$} \label{Final_2}
For all integer $k \geq 0$:
\begin{multline} \label{eq:C^z_2k+1_parc}
\sum_{j=1}^{n}\frac{1}{(a j+b)^{2k+1}}\cos{\frac{2\pi(a j+b)}{z}}=-\frac{1}{2b^{2k+1}}\left(\cos{\frac{2\pi b}{z}}-\sum_{j=0}^{k}\frac{(-1)^j}{(2j)!}\left(\frac{2\pi b}{z}\right)^{2j}\right)\\+\frac{1}{2(a n+b)^{2k+1}}\left(\cos{\frac{2\pi(a n+b)}{z}}-\sum_{j=0}^{k}\frac{(-1)^j}{(2j)!}\left(\frac{2\pi(a n+b)}{z}\right)^{2j}\right)
\\+\sum_{j=0}^{k}\frac{(-1)^{k-j}}{(2k-2j)!}\left(\frac{2\pi}{z}\right)^{2k-2j}\HP_{2j+1}(n)\\+\frac{(-1)^k}{2(2k)!}\left(\frac{2\pi}{z}\right)^{2k+1}\int_{0}^{1}(1-u)^{2k}\left(\cos{\frac{2\pi(a n+b)u}{z}}-\cos{\frac{2\pi b u}{z}}\right)\cot{\frac{\pi a u}{z}}\,du 
\end{multline}

For all integer $k \geq 1$:
\begin{multline} \label{eq:S^z_2k_parc}
\sum_{j=1}^{n}\frac{1}{(a j+b)^{2k}}\sin{\frac{2\pi(a j+b)}{z}}=-\frac{1}{2b^{2k}}\left(\sin{\frac{2\pi b}{z}}-\sum_{j=0}^{k-1}\frac{(-1)^j}{(2j+1)!}\left(\frac{2\pi b}{z}\right)^{2j+1}\right)\\+\frac{1}{2(a n+b)^{2k}}\left(\sin{\frac{2\pi(a n+b)}{z}}-\sum_{j=0}^{k-1}\frac{(-1)^j}{(2j+1)!}\left(\frac{2\pi(a n+b)}{z}\right)^{2j+1}\right)
\\-\sum_{j=0}^{k-1}\frac{(-1)^{k-j}}{(2k-1-2j)!}\left(\frac{2\pi}{z}\right)^{2k-1-2j}\HP_{2j+1}(n)
\\-\frac{(-1)^k}{2(2k-1)!}\left(\frac{2\pi}{z}\right)^{2k}\int_{0}^{1}(1-u)^{2k-1}\left(\cos{\frac{2\pi(a n+b)u}{z}}-\cos{\frac{2\pi b u}{z}}\right)\cot{\frac{\pi a u}{z}}\,du 
\end{multline} 

\subsubsection{The limits of $C^z_{2k+1}(n)$ and $S^z_{2k}(n)$} \label{dificil}
The limits of $C^z_{2k+1}(n)$ and $S^z_{2k}(n)$ for real $\abs{z}\ge 1$ are $C^z_{2k+1}$ and $S^z_{2k}$, which are given, respectively, by:
\begin{multline} \label{eq:C^z_{2k+1}}
\sum_{j=1}^{\infty}\frac{1}{j^{2k+1}}\cos{\frac{2\pi j}{z}}=\sum_{j=1}^{k}\frac{(-1)^{k-j}}{(2k-2j)!}\left(\frac{2\pi}{z}\right)^{2k-2j}\zeta(2j+1)+\frac{(-1)^k}{(2k)!}\left(\frac{2\pi}{z}\right)^{2k}\log{\abs{z}}\\-\frac{(-1)^k}{2(2k)!}\left(\frac{2\pi}{z}\right)^{2k+1}\int_{0}^{1}(1-u)^{2k}\cot{\frac{\pi u}{z}}-z(1-u)\cot{\pi u}\,du 
\end{multline}
\begin{multline} \label{eq:S^z_{2k}}
\sum_{j=1}^{\infty}\frac{1}{j^{2k}}\sin{\frac{2\pi j}{z}}=-\sum_{j=1}^{k-1}\frac{(-1)^{k-j}}{(2k-1-2j)!}\left(\frac{2\pi}{z}\right)^{2k-1-2j}\zeta(2j+1)-\frac{(-1)^k }{(2k-1)!}\left(\frac{2\pi}{z}\right)^{2k-1}\log{\abs{z}}\\+\frac{(-1)^k}{2(2k-1)!}\left(\frac{2\pi}{z}\right)^{2k}\int_{0}^{1}(1-u)^{2k-1}\cot{\frac{\pi u}{z}}-z(1-u)\cot{\pi u}\,du 
\end{multline}
\indent The exception is $C^1_1=\infty$, since integral $\int_{0}^{1}\cot{\pi u}-(1-u)\cot{\pi u}\,du$ diverges, which implies that $H(n)$ diverges. These results are probably original.

\section{The limits of the integrals} \label{int_lim}

In \citena{GHNR}, the following theorems were introduced, whose validity is now fully extended. For all real $k \geq 0$ and real $z$:
\begin{equation} \nonumber
\textbf{Theorem 1}\lim_{n\to\infty}\int_{0}^{1}(1-u)^{k}\sin{\frac{2\pi n u}{z}}\cot{\frac{\pi u}{z}}\,du=
\begin{cases}
      1, & \text{if}\ k=0 \text{ and }\abs{z}=1\\
      \abs{z/2}, & \text{if }\abs{z}\ge 1
\end{cases} 
\end{equation}
\indent Another result needed is stated in the following theorem, which holds for all real $k \geq 0$ and real $\abs{z}\geq 1$ (except $k=0$ and $\abs{z}=1$, for which the integral does not converge):
\begin{equation} \nonumber
\textbf{Theorem 2}
\lim_{n\to\infty}\int_{0}^{1}(1-u)^{k}\cos{\frac{2\pi n u}{z}}\cot{\frac{\pi u}{z}}-z\,(1-u)\cos{2\pi n u}\cot{\pi u}\,du=\frac{z\log{\abs{z}}}{\pi} 
\end{equation}
\indent A direct consequence of Theorem 2 and of the various possible formulae for $H(n)$ (see \citena{GHNR}), the result below is useful to handle the half-integers in \secrefe{C_lim} and \secrefe{Final}:
\begin{equation} \nonumber
\lim_{n\to\infty}\int_{0}^{1}(1-u)^{k}\cos{\frac{2\pi n u}{z}}\cot{\frac{\pi u}{z}}-\frac{z}{2}\,(1-u)\cos{\pi n u}\cot{\frac{\pi u}{2}}\,du=\frac{z}{\pi}\log{\frac{\abs{z}}{2}} 
\end{equation}
\indent Note that, at the time this paper was first released, these results were conjectures and had not been formally demonstrated yet.\\

Though these limits should not converge for non-real complex $z$, when they are linearly combined like $l_1+\ii\,l_2$ their infinities cancel out giving a finite value. This property is what allows the final formula from section \secrefe{Final} to converge nearly always, even when the parameters are not real.

\section{$\HP(n)$ asymptotic behavior}
Herein a relation between $\HP(n)$ and $H(n)$ is established.\\ 

For this exercise, the sine-based $\HP(n)$ formula from \citena{CHP} is used:
\begin{equation} \nonumber
\sum_{j=1}^{n}\frac{1}{j+b}=-\frac{1}{2b}+\frac{1}{2(n+b)}+\frac{\pi}{\sin{2\pi b}}\int_{0}^{1}\left(\sin{2\pi(n+b)u}-\sin{2\pi b u}\right)\cot{\pi u}\,du 
\end{equation}
\indent The first sine in the integral can be expanded (with the sine addition formula, $\sin(x+y)=\sin{x}\cos{y}+\cos{x}\sin{y}$), giving:
\begin{multline} \nonumber
\sum_{j=1}^{n}\frac{1}{j+b}=-\frac{1}{2b}+\frac{1}{2(n+b)}\\+\frac{\pi}{\sin{2\pi b}}\int_{0}^{1}\left(\cos{2\pi bu}\sin{2\pi n u}+\sin{2\pi bu}\cos{2\pi n u}-\sin{2\pi b u}\right)\cot{\pi u}\,du  
\end{multline}
\indent By taking the first term of the function in the integral, making a change of variables ($u$ to $1-u$), and expanding $\cos{2\pi b(1-u)}$ (with the cosine addition formula, $\cos(x+y)=\cos{x}\cos{y}-\sin{x}\sin{y}$), it follows from Theorem 1 that:
\begin{multline} \label{eq:1o_lim}
\lim_{n\to\infty}\frac{\pi}{\sin{2\pi b}}\int_{0}^{1}\cos{2\pi b(1-u)}\sin{2\pi n(1-u)}\cot{\pi(1-u)}\,du \\=\frac{\pi}{2}\left(\cot{2\pi b}+\csc{2\pi b}\right)=\frac{\pi}{2}\cot{\pi b}
\end{multline}
\indent The same steps are repeated for the remaining terms of the function in the integral. A change of variables is made and $\sin{2\pi b(1-u)}$ is expanded with the sine addition formula. But this time, when using Theorem 2 one must exclude the case $k=0$ and $z=1$ (since that integral does not converge), which gives:
\begin{multline} \label{eq:1o_lim_2}
\int_{0}^{1}\left(\sin{2\pi b(1-u)}\cos{2\pi n(1-u)}-\sin{2\pi b(1-u)}\right)\cot{\pi(1-u)}\,du=\\
\sin{2\pi b}\int_{0}^{1}(\cos{2\pi bu}-1-u(\cos{2\pi b}-1))\cos{2\pi n(1-u)}\cot{\pi(1-u)}\,du \\
-\cos{2\pi b}\int_{0}^{1}(\sin{2\pi bu}-u\sin{2\pi b})\cos{2\pi n(1-u)}\cot{\pi(1-u)}\,du \\
+\int_{0}^{1}(-\sin{2\pi b(1-u)}+(\sin{2\pi b})(1-u)\cos{2\pi n(1-u)})\cot{\pi(1-u)}\,du
\end{multline}
\indent The two first integrals on the right-hand side cancel out, per Theorem 2, when $n$ goes to infinity, leaving only the third integral to evaluate.\\

But looking back at the expression for $H(n)$ from \citena{GHNR}, one notices that it matches part of the last integral:
\begin{equation} \label{eq:H(n)}
\sum_{j=1}^{n}\frac{1}{j}=\frac{1}{2n}+\pi\int_{0}^{1} u\left(1-\cos{2\pi n (1-u)}\right)\cot{\pi(1-u)}\,du \text{,}
\end{equation}
\noindent which means the last integral can be further split:
\begin{multline} \nonumber
\frac{\pi}{\sin{2\pi b}}\int_{0}^{1}(-\sin{2\pi b(1-u)}+(\sin{2\pi b})(1-u)\cos{2\pi n(1-u)})\cot{\pi(1-u)}\,du=\\
\frac{\pi}{\sin{2\pi b}}\int_{0}^{1}(-\sin{2\pi b(1-u)}-u\sin{2\pi b}+\sin{2\pi b}\cos{2\pi n(1-u)})\cot{\pi(1-u)}\,du\\
+\pi\int_{0}^{1} u\left(1-\cos{2\pi n (1-u)}\right)\cot{\pi(1-u)}\,du 
\end{multline}
\indent At this point, there is only the limit of the first integral on the right-hand side left to evaluate, but fortunately that integral is constant for all integer $n$\footnote{It stems from $\int_{0}^{1}(1-\cos{2\pi n\,u})\cot{\pi\,u}\,du=0$ for all integer $n$.}. Therefore, after simplifying \eqrefe{1o_lim} further, one concludes that for sufficiently large $n$:
\begin{equation} \label{eq:approx}
\sum_{j=1}^{n}\frac{1}{j+b} \sim -\frac{1}{2b}+\frac{\pi}{2}\cot{\pi b}-\pi\int_{0}^{1}\left(\frac{\sin{2\pi b u}}{\sin{2\pi b}}-u\right)\cot{\pi u}\,du+H(n) 
\end{equation}
\indent Coincidentally, the above integral is identical to the generating function of the zeta function at the odd integers, that was created in \citena{GHNR}:
\begin{equation} \nonumber
\sum_{k=1}^{\infty}\zeta(2k+1)x^{2k+1}=-\pi x\int_{0}^{1}\left(\frac{\sin{2\pi x u}}{\sin{2\pi x}}-u\right)\cot{\pi u}\,du 
\end{equation}
\indent That means that for sufficiently large $n$ and $0<\abs{b}<1$, one can write the interesting approximation:
\begin{equation} \nonumber
\sum_{j=1}^{n}\frac{1}{j+b} \sim H(n)-\sum_{k=2}^{\infty}(-1)^k\zeta(k)b^{k-1}=H(n)-H(b)
\end{equation}
\indent Now, since formula \eqrefe{approx} clearly does not hold at the half-integers, for such $b$ one can resort to a different integral representation of the generating function of $\zeta(2k+1)$\citesup{GHNR}, which leads to:
\begin{equation} \nonumber
\sum_{j=1}^{n}\frac{1}{j+b} \sim -\frac{1}{2b}+\frac{\pi}{2}\int_{0}^{1}\left(-1+u+\cos{\pi b u}\right)\cot{\frac{\pi u}{2}}\,du+H(n) \text{,}
\end{equation}
\noindent since $\cot{\pi b}$ is zero for all half-integer $b$.

\section{The full Fourier series} \label{new_lim}
Although the expressions of $C^z_{k}(a,b,n)$ and $S^z_{k}(a,b,n)$ hold for all positive integers $k$ and $n$ and complex $z$, $a$ and $b$, the limits found next are constrained by the requirements of the Theorems 1 and 2 from section \secrefe{int_lim}.\\

Without loss of generality, let $a=1$ to simplify the calculations:
\begin{equation} \nonumber
C^z_{k}(b,n)=\sum_{j=1}^{n}\frac{1}{(j+b)^k}\cos{\frac{2\pi(j+b)}{z}} \text{ and } S^z_{k}(b,n)=\sum_{j=1}^{n}\frac{1}{(j+b)^k}\sin{\frac{2\pi(j+b)}{z}}
\end{equation}
\indent And since $k=0$ and $\abs{z}=1$ leads to trivial cases, they are not accounted for in the following demonstration (so the final formulae may not be true for $k=0$ and $\abs{z}=1$).

\subsection{The limit of $C^z_{2k+1}(b,n)$} \label{C_lim}
The limit of $C^z_{2k+1}(b,n)$ is much harder to obtain than the limit of $S^z_{2k+1}(b,n)$, which should come as no surprise, given the limits that were shown in sections \secrefe{facil} and \secrefe{dificil}, for the particular cases.\\

When the limit of $C^z_{2k+1}(b,n)$, a particular case of formula  \eqrefe{C^z_2k+1_parc} with $a=1$, as $n$ goes to infinity is taken, the only challenging terms are the ones that diverge, that is, the $\HP(n)$ (that appears in the last summation) and the integral. The intermediate result is shown below:
\begin{multline} \nonumber
\lim_{n\to\infty}C^z_{2k+1}(b,n)=-\frac{1}{2b^{2k+1}}\left(\cos{\frac{2\pi b}{z}}-\sum_{j=0}^{k}\frac{(-1)^j }{(2j)!}\left(\frac{2\pi b}{z}\right)^{2j}\right)
\\+\sum_{j=1}^{k}\frac{(-1)^{k-j}}{(2k-2j)!}\left(\frac{2\pi}{z}\right)^{2k-2j}\zeta(2j+1,b+1)
\\+\lim_{n\to\infty}\frac{(-1)^k}{2(2k)!}\left(\frac{2\pi}{z}\right)^{2k+1}\left(\frac{z}{\pi}\sum_{j=1}^{n}\frac{1}{j+b}+\int_{0}^{1}(1-u)^{2k}\left(\cos{\frac{2\pi(n+b)u}{z}}-\cos{\frac{2\pi b u}{z}}\right)\cot{\frac{\pi u}{z}}\,du\right)
\end{multline}
\indent Now, if one recalls the approximation found for $\HP(n)$ in \eqrefe{approx}, $\HP(n) \sim H(n)+c$ for large $n$ (where $c$ is the term that does not depend on $n$). Therefore, one only needs to solve the limit:
\begin{multline} \nonumber
\lim_{n\to\infty}\left(\frac{z}{\pi}(c+H(n))+\int_{0}^{1}(1-u)^{2k}\left(\cos{\frac{2\pi(n+b)u}{z}}-\cos{\frac{2\pi b u}{z}}\right)\cot{\frac{\pi u}{z}}\,du\right)
\end{multline}
\indent The cosine in the integral can be expanded with the cosine addition formula giving:
\begin{equation} \nonumber
\int_{0}^{1}(1-u)^{2k}\left(\cos{\frac{2\pi bu}{z}}\left(-1+\cos{\frac{2\pi nu}{z}}\right)-\sin{\frac{2\pi b u}{z}}\sin{\frac{2\pi n u}{z}}\right)\cot{\frac{\pi u}{z}}\,du
\end{equation}
\indent But due to Theorem 1, the below limit is zero (to see that, one just needs to expand the first sine --- note $u$ was changed for $1-u$):
\begin{equation} \nonumber
\lim_{n\to\infty}\int_{0}^{1}u^{2k}\sin{\frac{2\pi b(1-u)}{z}}\sin{\frac{2\pi n(1-u)}{z}}\cot{\frac{\pi(1-u)}{z}}\,du=0
\end{equation}
\indent Now, by replacing $H(n)$ with its formula from equation \eqrefe{H(n)} and adding it to what is left in the integral, one has:
\begin{equation} \nonumber
\int_{0}^{1}(1-u)^{2k}\cos{\frac{2\pi bu}{z}}\left(-1+\cos{\frac{2\pi nu}{z}}\right)\cot{\frac{\pi u}{z}}+z(1-u)\left(1-\cos{2\pi n u}\right)\cot{\pi u}\,du
\end{equation}
\indent Looking at Theorem 2, one can recombine the terms conveniently into an integral that converges as $n$ goes to infinity:
\begin{equation} \nonumber
\lim_{n\to\infty}\int_{0}^{1}(1-u)^{2k}\cos{\frac{2\pi bu}{z}}\cos{\frac{2\pi nu}{z}}\cot{\frac{\pi u}{z}}-z(1-u)\cos{2\pi n u}\cot{\pi u}\,du=\frac{z\log{\abs{z}}}{\pi} \text{,}
\end{equation}
\noindent which is justified by the following (only one term shown):
\begin{multline} \nonumber
\cos{\frac{2\pi b}{z}}\int_{0}^{1}u^{2k}\cos{\frac{2\pi bu}{z}}\cos{\frac{2\pi n(1-u)}{z}}\cot{\frac{\pi(1-u)}{z}}\\-z\cos{\frac{2\pi b}{z}}\,u\cos{2\pi n(1-u)}\cot{\pi(1-u)}\,du \rightarrow \left(\cos{\frac{2\pi b}{z}}\right)^2\frac{z\log{\abs{z}}}{\pi} \text{,}
\end{multline}
\noindent whereas the remaining integral converges on its own.\\

Let us summarize the result. The below limit does not change if one picks different formulae for $H(n)$, which is useful to determine how the formula changes for the half-integers $b$:
\begin{multline} \nonumber
\lim_{n\to\infty}\left(\frac{z}{\pi}H(n)+\int_{0}^{1}(1-u)^{2k}\left(\cos{\frac{2\pi(n+b)u}{z}}-\cos{\frac{2\pi b u}{z}}\right)\cot{\frac{\pi u}{z}}\,du\right)\\
=\frac{z\log{\abs{z}}}{\pi}
-\int_{0}^{1}(1-u)^{2k}\cos{\frac{2\pi bu}{z}}\cot{\frac{\pi u}{z}}-z(1-u)\cot{\pi u}\,du\\=\frac{z}{\pi}\log{\frac{\abs{z}}{2}}-\int_{0}^{1}(1-u)^{2k}\cos{\frac{2\pi bu}{z}}\cot{\frac{\pi u}{z}}-\frac{z}{2}(1-u)\cot{\frac{\pi u}{2}}\,du 
\end{multline}

\subsubsection{Non-integer $2b$}
\indent After one puts everything together, the conclusion is that for all integer $k\ge 0$ and real $\abs{z}\ge 1$:
\begin{multline} \label{eq:final_non_int}
\sum_{j=1}^{\infty}\frac{1}{(j+b)^{2k+1}}\cos{\frac{2\pi(j+b)}{z}}=-\frac{1}{2b^{2k+1}}\left(\cos{\frac{2\pi b}{z}}-\sum_{j=0}^{k-1}\frac{(-1)^j}{(2j)!}\left(\frac{2\pi b}{z}\right)^{2j}\right)
\\+\sum_{j=1}^{k}\frac{(-1)^{k-j}}{(2k-2j)!}\left(\frac{2\pi}{z}\right)^{2k-2j}\zeta(2j+1,b+1)
+\frac{(-1)^k\pi}{2(2k)!}\left(\frac{2\pi}{z}\right)^{2k}\cot{\pi b}+\frac{(-1)^k}{(2k)!}\left(\frac{2\pi}{z}\right)^{2k}\log{\abs{z}}\\-\frac{(-1)^k}{2(2k)!}\left(\frac{2\pi}{z}\right)^{2k+1}\int_{0}^{1}(1-u)^{2k}\cos{\frac{2\pi bu}{z}}\cot{\frac{\pi u}{z}}+z\left(-1+\frac{\sin{2\pi b u}}{\sin{2\pi b}}\right)\cot{\pi u}\,du
\end{multline}
\indent As one can see, it takes a really convoluted function to generate this simple Fourier series.

\subsubsection{Half-integer $b$}
At the half-integers $b$, the formula reduces to:
\begin{multline} \label{eq:final_half_int}
\sum_{j=1}^{\infty}\frac{1}{(j+b)^{2k+1}}\cos{\frac{2\pi(j+b)}{z}}=-\frac{1}{2b^{2k+1}}\left(\cos{\frac{2\pi b}{z}}-\sum_{j=0}^{k-1}\frac{(-1)^j}{(2j)!}\left(\frac{2\pi b}{z}\right)^{2j}\right)
\\+\sum_{j=1}^{k}\frac{(-1)^{k-j}}{(2k-2j)!}\left(\frac{2\pi}{z}\right)^{2k-2j}\zeta(2j+1,b+1)+\frac{(-1)^k}{(2k)!}\left(\frac{2\pi}{z}\right)^{2k}\log{\frac{\abs{z}}{2}}\\-\frac{(-1)^k}{2(2k)!}\left(\frac{2\pi}{z}\right)^{2k+1}\int_{0}^{1}(1-u)^{2k}\cos{\frac{2\pi bu}{z}}\cot{\frac{\pi u}{z}}-\frac{z}{2}\cos{\pi b u}\cot{\frac{\pi u}{2}}\,du
\end{multline}

\subsubsection{Integer $b$}
For integer $b$:
\begin{multline} \label{eq:final_int}
\sum_{j=1}^{\infty}\frac{1}{(j+b)^{2k+1}}\cos{\frac{2\pi(j+b)}{z}}=-\frac{1}{2b^{2k+1}}\left(\cos{\frac{2\pi b}{z}}-\sum_{j=0}^{k}\frac{(-1)^j}{(2j)!}\left(\frac{2\pi b}{z}\right)^{2j}\right)
\\+\sum_{j=1}^{k}\frac{(-1)^{k-j}}{(2k-2j)!}\left(\frac{2\pi}{z}\right)^{2k-2j}\zeta(2j+1,b+1)-\frac{(-1)^k}{(2k)!}\left(\frac{2\pi}{z}\right)^{2k}H(b)+\frac{(-1)^k}{(2k)!}\left(\frac{2\pi}{z}\right)^{2k}\log{\abs{z}}\\-\frac{(-1)^k}{2(2k)!}\left(\frac{2\pi}{z}\right)^{2k+1}\int_{0}^{1}(1-u)^{2k}\cos{\frac{2\pi bu}{z}}\cot{\frac{\pi u}{z}}-z(1-u)\cot{\pi u}\,du
\end{multline}

\subsection{The limit of $S^z_{2k+1}(b,n)$}
In the case of $S^z_{2k+1}(b,n)$, regardless of integer or half-integer, one has:
\begin{multline} \nonumber
\lim_{n\to\infty}S^z_{2k+1}(b,n)=-\frac{1}{2b^{2k+1}}\left(\sin{\frac{2\pi b}{z}}-\sum_{j=0}^{k-1}\frac{(-1)^j}{(2j+1)!}\left(\frac{2\pi b}{z}\right)^{2j+1}\right)
\\+\sum_{j=1}^{k}\frac{(-1)^{k-j}}{(2k+1-2j)!}\left(\frac{2\pi}{z}\right)^{2k+1-2j}\zeta(2j,b+1)\\
+\lim_{n\to\infty}\frac{(-1)^k}{2(2k)!}\left(\frac{2\pi}{z}\right)^{2k+1}\int_{0}^{1}(1-u)^{2k}\left(\sin{\frac{2\pi(n+b)u}{z}}-\sin{\frac{2\pi b u}{z}}\right)\cot{\frac{\pi u}{z}}\,du
\end{multline}
\indent This one is much simpler and the limit of the integral can be easily deduced by means of the Theorem 1, without even having to expand the sine in the integral (with the sine addition formula).\\

Thus, for all integer $k\ge 0$ and real $\abs{z}\ge 1$:
\begin{multline} \label{eq:S^z_2k+1_final}
\sum_{j=1}^{\infty}\frac{1}{(j+b)^{2k+1}}\sin{\frac{2\pi(j+b)}{z}}=-\frac{1}{2b^{2k+1}}\left(\sin{\frac{2\pi b}{z}}-\sum_{j=0}^{k-1}\frac{(-1)^j}{(2j+1)!}\left(\frac{2\pi b}{z}\right)^{2j+1}\right)
\\+\sum_{j=1}^{k}\frac{(-1)^{k-j}}{(2k+1-2j)!}\left(\frac{2\pi}{z}\right)^{2k+1-2j}\zeta(2j,b+1)+\frac{(-1)^k\abs{z}}{4(2k)!}\left(\frac{2\pi}{z}\right)^{2k+1}\\-\frac{(-1)^k}{2(2k)!}\left(\frac{2\pi}{z}\right)^{2k+1}\int_{0}^{1}(1-u)^{2k}\sin{\frac{2\pi b u}{z}}\cot{\frac{\pi u}{z}}\,du
\end{multline}

\subsection{$C^z_{2k}(b)$ and $S^z_{2k}(b)$}

The next two formulae, $C^z_{2k}(b)$ and $S^z_{2k}(b)$, are analogs and do not require further explanations:
\begin{multline} \label{eq:C^z_2k(b)}
\sum_{j=1}^{\infty}\frac{1}{(j+b)^{2k}}\cos{\frac{2\pi(j+b)}{z}}=-\frac{1}{2b^{2k}}\left(\cos{\frac{2\pi b}{z}}-\sum_{j=0}^{k-1}\frac{(-1)^j}{(2j)!}\left(\frac{2\pi b}{z}\right)^{2j}\right)
\\+\sum_{j=1}^{k}\frac{(-1)^{k-j}}{(2k-2j)!}\left(\frac{2\pi}{z}\right)^{2k-2j}\zeta(2j,b+1)+\frac{(-1)^k\abs{z}}{4(2k-1)!}\left(\frac{2\pi}{z}\right)^{2k}\\-\frac{(-1)^k}{2(2k-1)!}\left(\frac{2\pi}{z}\right)^{2k}\int_{0}^{1}(1-u)^{2k-1}\sin{\frac{2\pi b u}{z}}\cot{\frac{\pi u}{z}}\,du
\end{multline}

\begin{multline} \nonumber
\sum_{j=1}^{\infty}\frac{1}{(j+b)^{2k}}\sin{\frac{2\pi(j+b)}{z}}=-\frac{1}{2b^{2k}}\left(\sin{\frac{2\pi b}{z}}-\sum_{j=0}^{k-2}\frac{(-1)^j}{(2j+1)!}\left(\frac{2\pi b}{z}\right)^{2j+1}\right)
\\-\sum_{j=1}^{k-1}\frac{(-1)^{k-j}}{(2k-1-2j)!}\left(\frac{2\pi}{z}\right)^{2k-1-2j}\zeta(2j+1,b+1)
\\-\frac{(-1)^k\pi}{2(2k-1)!}\left(\frac{2\pi}{z}\right)^{2k-1}\cot{\pi b}-\frac{(-1)^k\log{\abs{z}}}{(2k-1)!}\left(\frac{2\pi}{z}\right)^{2k-1}\\+\frac{(-1)^k}{2(2k-1)!}\left(\frac{2\pi}{z}\right)^{2k}\int_{0}^{1}(1-u)^{2k-1}\cos{\frac{2\pi bu}{z}}\cot{\frac{\pi u}{z}}+z\left(-1+\frac{\sin{2\pi b u}}{\sin{2\pi b}}\right)\cot{\pi u}\,du
\end{multline}

\section{Lerch's $\Phi$ at the positive integers}
Now a proper formula for the Lerch transcendent function, $\Phi(e^z,k,b)$, at the positive integers $k$ can be created. To avoid repetition, note that the expressions obtained for the limits of the partial formulae hold for all complex $z$, except $z$ such that both $\Re{(z)}>=0$ and $\abs{\Im{(z)}}>2\pi$.

\subsection{Partial Lerch $\Phi$ sums, $E^{z}_{k}(b,n)$}
It is straightforward to derive an expression for the partial sums of the Lerch $\Phi$ function using the formulae from \secrefe{Final_1} and \secrefe{Final_2}. If $\ii$ is the imaginary unit, one just makes:
\begin{equation} \nonumber
E^{2\pi\ii/z}_{k}(b,n)=\sum_{j=1}^{n}\frac{e^{2\pi\ii(j+b)/z}}{(j+b)^{k}}=C^z_{k}(b,n)+\ii\,S^z_{k}(b,n) 
\end{equation}
\indent By omitting the calculations and making a simple transformation ($z:=2\pi\ii/z$) (to bring the variables into the domain of the real numbers), one can produce a single formula for both the odd and even powers:
\begin{multline} \nonumber 
\sum _{j=1}^{n}\frac{e^{z(j+b)}}{(j+b)^{k}}=-\frac{e^{z b}}{2b^{k}}+\frac{e^{z(n+b)}}{2(n+b)^{k}}+\frac{1}{2b^{k}}\sum_{j=0}^{k}\frac{(z b)^j}{j!}-\frac{1}{2(n+b)^{k}}\sum_{j=0}^{k}\frac{(z (n+b))^j}{j!}\\+\sum_{j=1}^{k}\frac{z^{k-j}}{(k-j)!}\HP_{j}(n)+\frac{z^{k}}{2(k-1)!}\int_0^1(1-u)^{k-1}\left(e^{z(n+b) u}-e^{z b u}\right)\coth{\frac{z u}{2}}\,du 
\end{multline}
\indent From this new equation, it is easy to see that as $n$ goes to infinity, the summation on the left-hand side converges only if $\Re{(z)}<0$. However, one can obtain an analytic continuation for this summation, by removing the second term on the right-hand side, which diverges if $\Re{(z)}>0$. Perhaps not surprisingly, this analytic continuation coincides with the Lerch $\Phi$ function.

\subsection{Partial polylogarithm sums, $E^{z}_{k}(0,n)$} \label{part}

When $b=0$, one has an interesting particular case:
\begin{multline} \label{eq:ab}
\sum _{j=1}^n \frac{e^{z j}}{j^{k}}=\frac{e^{z n}}{2n^{k}}-\frac{1}{2n^{k}}\sum_{j=0}^{k}\frac{(z n)^j}{j!}+\sum_{j=1}^{k}\frac{z^{k-j}}{(k-j)!}H_{j}(n)\\
+\frac{z^{k}}{2(k-1)!}\int_0^1(1-u)^{k-1}\left(e^{z n u}-1\right)\coth{\frac{z u}{2}}\,du 
\end{multline}\\
\indent To obtain this expression, the limit as $b$ tends to zero is taken:
\begin{equation} \nonumber 
\lim_{b\to 0}-\frac{e^{z b}}{2b^{k}}+\frac{1}{2b^{k}}\sum_{j=0}^{k}\frac{(z b)^j}{j!}=0 
\end{equation}

\subsection{Lerch's $\Phi$} \label{Final}
The limits found in section \secrefe{new_lim} allow to create a formula for the infinite series below:
\begin{equation} \nonumber
\sum_{j=1}^{\infty}\frac{e^{2\pi\ii (j+b)/z}}{(j+b)^{k}}=\lim_{n\to\infty}C^z_{k}(b,n)+\ii S^z_{k}(b,n)
\end{equation}

\subsubsection{Non-integer $2b$}
After all the calculations are carried out, one finds that for all integer $k \ge 1$:
\begin{multline} \label{eq:Lerch_final_2b} 
\sum _{j=1}^{\infty}\frac{e^{z(j+b)}}{(j+b)^{k}}=-\frac{1}{2b^{k}}\left(e^{z b}-\sum_{j=0}^{k-2}\frac{(z b)^j}{j!}\right)+\sum_{j=2}^{k}\frac{z^{k-j}}{(k-j)!}\zeta(j,b+1)\\+\frac{\pi\,z^k}{2(k-1)!}\cot{\pi b}-\frac{z^{k-1}}{(k-1)!}\log{\left(-\frac{z}{2\pi}\right)}\\-\frac{z^{k}}{2(k-1)!}\int_0^1(1-u)^{k-1}e^{z\,b\,u}\coth{\frac{z u}{2}}+\frac{2\pi}{z}\left(-1+\frac{\sin{2\pi b u}}{\sin{2\pi b}}\right)\cot{\pi u}\,du 
\end{multline}
\indent The infinite series on the left-hand side converges whenever $\Re{(z)}<0$, whereas the expression on the right-hand side, $E^{z}_{k}(b)$, is well defined always, except when $2b$ is an integer or $z$ lies outside the domain. At $z=0$, although improper, the expression has a limit.\\

This series is related to the actual Lerch $\Phi$ function by a very simple relation:
\begin{equation} \nonumber
\Phi(e^z,k,b)=\frac{1}{b^k}+e^{-z\,b}\sum _{j=1}^{\infty}\frac{e^{z(j+b)}}{(j+b)^{k}} \text{}
\end{equation}

\subsubsection{Half-integer $b$}
For half-integer $b$, the formula is only slightly different. For all integer $k \ge 1$:
\begin{multline}  \label{eq:Lerch_final_half_int} 
\sum _{j=1}^{\infty}\frac{e^{z(j+b)}}{(j+b)^{k}}=-\frac{1}{2b^{k}}\left(e^{z b}-\sum_{j=0}^{k-2}\frac{(z b)^j}{j!}\right)+\sum_{j=2}^{k}\frac{z^{k-j}}{(k-j)!}\zeta(j,b+1)\\-\frac{z^{k-1}}{(k-1)!}\log{\left(-\frac{z}{\pi}\right)}-\frac{z^{k}}{2(k-1)!}\int_0^1(1-u)^{k-1}e^{z\,b\,u}\coth{\frac{z u}{2}}-\frac{\pi}{z}\cos{\pi b u}\cot{\frac{\pi u}{2}}\,du 
\end{multline}

\subsubsection{Integer $b$}
When $b$ is a positive integer, $E^{z}_{k}(b)$ becomes an incomplete polylog series, which is covered next. Hence, it is very simple to derive its formula, one just needs to subtract the missing term from the full polylog. A similar reasoning is used if $b$ is a negative integer.\\

Nonetheless, the formula when $b$ is a positive integer is:
\begin{multline}  \label{eq:Lerch_final_int}
\sum _{j=1}^{\infty}\frac{e^{z(j+b)}}{(j+b)^{k}}=-\frac{1}{2b^{k}}\left(e^{z b}-\sum_{j=0}^{k-2}\frac{(z b)^j}{j!}\right)+\sum_{j=2}^{k}\frac{z^{k-j}}{(k-j)!}\zeta(j,b+1)\\-\frac{z^{k-1}}{(k-1)!}\log{\left(-\frac{z}{2\pi}\right)}-\frac{z^{k-1}}{(k-1)!}\left(H(b)-\frac{1}{2b}\right)\\-\frac{z^{k}}{2(k-1)!}\int_0^1(1-u)^{k-1}e^{z\,b\,u}\coth{\frac{z u}{2}}-\frac{2\pi}{z}(1-u)\cot{\pi u}\,du 
\end{multline}

\subsection{The polylogarithm, $\mathrm{Li}_{k}(e^{z})$}
The limit of $E^{z}_{k}(0,n)$ when $n$ tends to infinity is the limit of the expression that was just found when $b$ tends to zero, and it relies on the following two notable limits:
\begin{equation} \nonumber
\lim_{b\to 0}-\frac{1}{2b^{k}}\left(e^{z b}-\sum_{j=0}^{k-2}\frac{(z b)^j}{j!}\right)+\frac{\pi\,z^k}{2(k-1)!}\cot{\pi b}=-\frac{z^k}{2\,k!} \text{, and } \lim_{b\to 0}\frac{\sin{2\pi b u}}{\sin{2\pi b}}=u
\end{equation}
\indent Therefore, for all integer $k\ge 1$:
\begin{multline} \nonumber
\sum _{j=1}^{\infty} \frac{e^{z j}}{j^{k}}=-\frac{z^{k-1}}{(k-1)!}\log{\left(-\frac{z}{2\pi}\right)}+\sum_{\substack{j=0 \\ j\neq 1}}^{k}\frac{z^{k-j}}{(k-j)!}\zeta(j)\\
-\frac{z^{k}}{2(k-1)!}\int_{0}^{1}(1-u)^{k-1}\coth{\frac{z u}{2}}-\frac{2\pi}{z}(1-u)\cot{\pi u}\,du
\end{multline}
\indent This infinite series is known as the polylogarithm, $\mathrm{Li}_{k}(e^{z})$, and the formula on the right-hand side provides an analytic continuation for it for when $\Re{(z)}>0$ (as long as $\Im{(z)} \le 2\pi$).\\

Note how the first limit fits perfectly into the summation (together with the other $\zeta(j)$ values, except for the singularity). It is easy to show that when $z$ goes to zero, the formula found goes to $\zeta(k)$, if $k\ge 2$.\\

Besides, this $\mathrm{Li}_{k}(e^{z})$ formula allows to deduce the following power series for $e^z$:
\begin{equation} \nonumber
\lim_{k\to\infty} \sum_{\substack{j=2}}^{k}\frac{z^{k-j}}{(k-j)!}\zeta(j)=e^z
\end{equation}

\subsection{The Hurwitz zeta function, $\zeta(-k,b)$}
From the literature, it is well known that the Hurwitz zeta function is related to the polylog function by means of a relatively simple relation:
\begin{equation} \nonumber
\frac{(2\pi)^k}{\Gamma(k)}\zeta(1-k,b)=\ii^{-k}\,\mathrm{Li}_{k}(e^{2\pi\ii\,b})+\ii^{k}\,\mathrm{Li}_{k}(e^{-2\pi\ii\,b}) \text{,}
\end{equation}
\noindent which holds if $\abs{\Re{(b)}}\le 1$.\\

Since the formula found for $\mathrm{Li}_{k}(e^{z})$ holds at the positive integers $k$, this relation allows to obtain a formula for $\zeta(-k,b)$ that holds at the negative integers $-k$. Without showing the simple but long calculations involved, one concludes that despite the constraints of the initial relation, the below formula holds for every $b$:
\begin{equation} \nonumber
\zeta(-k,b)=\frac{b^k}{2}+2\,k!\,b^{k+1}\sum_{j=0}^{\floor{(k+1)/2}}\frac{(-1)^{j}(2\pi\,b)^{-2j}\zeta(2j)}{(k+1-2j)!}=-\frac{B_{k+1}(b)}{k+1} \text{,}
\end{equation}
\noindent where $B_{k+1}(b)$ are the Bernoulli polynomials, whose relation with the Hurwitz zeta is also known from literature. So that is an alternative expression for the Bernoulli polynomials.

\end{document}